\documentclass[12pt]{article}
\usepackage[cp1251]{inputenc}
\usepackage[T2A]{fontenc}
\usepackage[english]{babel}

\usepackage{amsmath,amssymb,amsfonts,amsthm}

\usepackage{graphicx,color}
\usepackage{xcolor}
\colorlet{blu1}{blue!70!black}
\colorlet{red1}{red!80}

\begin{document}



\begin{center}

{\large {{\bf Analogues of Jacobi and Weyl Theorems for Infinite-Dimensional Tori\footnote{This research was founded by Russian Scientific Foundation under grant 19-11-00320.}
}}}
\end{center}

\begin{center}
V.Zh. Sakbaev,\footnote{Steklov Mathematical Institute of Russian Academy of Sciences, fumi2003@mail.ru}\ I.V. Volovich\footnote{Steklov Mathematical Institute of Russian Academy of Sciences, volovich@mi-ras.ru}
\end{center}

{\centerline {\bf Abstract}}
Generalizations of the Jacobi and Weyl theorems on finite-dimensional linear flows to the case of linear flows on infinite-dimensional tori are presented.
Conditions  for periodicity, non-wandering,
ergodicity and transitivity of trajectories of an infinite-dimensional linear flow are obtained. It is shown that for infinite-dimensional linear flows there is a new type of trajectories that is absent in the finite-dimensional case.

 \noindent{\bf Keywords:}  Linear flow on infinite-dimensional tori; periodic trajectory; transitive trajectory; ergodic measure of a linear flow; nonwandering points.

\medskip
\noindent
{\bf Mathematical Subject Classification:}
28D05, 37A05, 37K10, 37N20

\section{Introduction}
{
An important class 
of infinite dimensional Hamiltonian systems
are infinite systems of oscillators.
In particular, every quantum system is unitarily equivalent to a system of classical harmonic oscillators \cite{Volovich21}.
}
The Koopman-von Neumann theory states that if the phase space of an infinite-dimensional Hamiltonian system is equipped with an invariant {measure} with respect to the Hamiltonian flow, then {this flow induces a group of unitary operators on the space of square integrable functions with respect to the given invariant measure \cite{Khrennikov, KS-18, Volovich19, S23}.
}
For this reason, many nonlinear PDEs admitting a Hamiltonian structure have been studied by introducing a measure invariant under the Hamiltonian flow \cite{Bourgain, Deng, Ponce, Sy}.
An ergodic theorem for flows in an infinite-dimensional Hilbert space was obtained in the work \cite{IvN}.

As proven in the paper \cite{Volovich21}, each unitary group {on} a Hilbert space
can be represented as {the} Hamiltonian flow of a system of harmonic oscillators in a {realization} of the Hilbert space equipped with a shift-invariant ({i.e. an invariant with respect to shift on any vector of the phase space})  symplectic form.
{Here the term realization of a complex Hilbert space means the introducing of the real structure on this space \cite{Khrennikov, KS-18}.}

The flow of a Hamiltonian system of oscillators corresponds to quantum dynamics in Hilbert space and has a complete family of first integrals in involution.
The invariant manifolds of level surfaces of the complete family of first integrals are invariant tori.

Properties (such as the existence of a periodic trajectory, {an everywhere dense (or topologically transitive) trajectory}) of the flow of a Hamiltonian system of oscillators (linear flow) on tori in a finite-dimensional symplectic space are studied in many papers \cite{Che, Dumas, 21, KH}. The existence of a complete system of conserved functionals for the periodic sine-Gordon equation and the corresponding infinite-dimensional tori is studied in \cite{Schw}. The development of the infinite-dimensional KAM theory is presented in  \cite{Che, KAM, Kuksin}.
A description of the measure on infinite-dimensional tori, together with the analysis of Fourier series, was carried out in \cite{Fufaev, Platonov}.

Weyl's theorem gives a criterion for the ergodicity of a linear flow on a finite-dimensional torus \cite{Che}. The ergodic theorem for infinite-dimensional Hilbert space was proved by I. von Neumann \cite{IvN}. The ergodicity of subgroups of orthogonal group in an infinite-dimensional Hilbert space is considered in \cite{Segal}.
{Generalizations of Weyl's theorem on ergodicity for the infinite-dimensional case were considered in \cite{KVV21}.}
An ergodic theorem on infinite-dimensional linear flow on an invariant torus equipped with the Kolmogorov measure is presented.

{In this paper we study the special case of a countable system of oscillators}. For a phase trajectory in a Hilbert space, criteria for periodicity (Theorem 2) and density everywhere on the invariant torus (Theorem 3) are obtained.
These results are a generalization to the infinite-dimensional case of Jacobi's theorem \cite{AA}, p. 114 (A1), which provides criteria for the periodicity of the system trajectory
$N$ harmonic oscillators and criteria for their density everywhere in the invariant torus. 
A generalization of the ergodicity criterion with respect to the Kolmogorov measure of a linear flow on an infinite-dimensional torus is presented (Theorem 4).

The system of $N$ harmonic oscillators is the Hamiltonian system in $2N$-dimensional Euclidean space $E={\mathbb R}^{2N}$ equipped with a shift-invariant symplectic form.
The Hamilton function of this Hamiltonian system is given by the following equality in the corresponding symplectic coordinates
\begin{equation}\label{HN}
H(q,p)=\sum\limits_{j=1}^N{{\lambda _j}\over 2}(q_j^2+p_j^2), \  (q,p)\in E,
\end{equation}
where $\{ \lambda _1,...,\lambda _N\}$ are
real numbers.
The system of harmonic oscillator (\ref{HN}) has the collection of $N$ first integrals  $I_j(q,p)=(q_j^2+p_j^2),\ j=1,...,N$.  For a given collection of $N$ positive numbers $r=\{ r_1,...,r_N\}$ the surface of first integral levels $T_r=\{ (q,p)\in E:\ I_j(q,p)=r_j,\ j=1,...,N\}$, is the smooth $N$-dimensional manifold which is  invariant with respect to the Hamiltonian flow (this manifold is called an invariant torus of the system of harmonic oscillators (\ref{HN})). The flow of a Harmonic oscillator system on its invariant torus is called linear flow \cite{KH}.

{\bf Definition 1.} (\cite{AA, KSF, KH}). {\it
A set of numbers $A\subset\mathbb R$ is said to be rationally commensurate if one of them is a linear combination of a finite set of others with rational coefficients.}

{\bf  Definition 2.}
{\it A set of numbers $A\subset \mathbb R$  is said to be strongly rationally commensurate if for any finite subset
$\{ a_1,...,a_m\}\subset A,\ m\geq 2,$ each of these numbers is a linear combination of $m-1$  other numbers with rational coefficients.}

{\bf Remark.}
The condition for rational commensurability of a finite set of numbers $A$ is equivalent to the condition for the resonance of a set of numbers $A$ (see \cite{21}, page 4).
Defining the resonance of finite multiplicity of a set of numbers $A$ in \cite{21} makes it possible to present conditions for strong commensurability of a finite set of numbers $A$.
A set of $N$ numbers $A$ is strongly commensurate if and only if it is a resonant set of maximum multiplicity $N-1$.

The following statement is known the Jacobi theorem for a finite system of oscillators (see \cite{AA}, p. 114).

{\bf Theorem (Jacobi)}.
{\it 1) A trajectory of the Hamiltonian system  (\ref{HN}) is periodic if and only if the set of frequencies $\{ \lambda _j,\ j=1,...,N\}$ is strongly rationally commensurable.

2) A trajectory of the Hamiltonian system (\ref{HN}) is everywhere dense on a $2N$-dimensional invariant torus $T=\{ I_k=c_k,\ k=1,...,N\}$ if and only if the collection of numbers does not satisfy the condition of rational commensurability,
i.e.
\begin{equation}\label{10}
\{\sum\limits_{k=1}^N\lambda _{k}n_k=0,\quad n_k\in {\mathbb Z}\} \quad \Rightarrow \quad \{ n_k=0,\ j=1,...,N\}.
\end{equation}
}

In this paper we prove the infinite-dimensional analogues of the Jacobi theorem (\cite{AA}), which provide criteria for the periodicity and  transitivity of a trajectory on an invariant torus of an infinite system of oscillators.
Criteria for the periodicity and transitivity of the trajectory of a countable system of oscillators are obtained (Theorems 2 and 3).

The Weyl theorem states that the linear flow on a $d$-dimensional torus
is ergodic with respect to Lebesgue measure on the invariant torus
if and only if the collection of frequencies is rationally independent \cite{Che}. We obtain an extension of this statement to the case of infinite-dimensional invariant torus in Theorem 4.

The Kolmogorov measure is studied on an invariant torus, which is a countable Cartesian product of circles. The Kolmogorov measure is defined as a countable product of probability measures on one-dimensional circles with a uniform distribution. It is proved (in Theorem 4) that the ergodicity with respect to the flow of a Hamiltonian system (\ref{HN}) of the Kolmogorov measure on a countable product of circles is equivalent to the transitivity of a linear flow on an invariant torus.

The naive direct reformulation of the Jacobi theorem for a countable system of oscillators has the following form "A trajectory of a countable system of oscillator is periodic if and only if the collection of oscillator frequencies is strongly rationally commensurable".
But in the infinite-dimensional case this is not true, as shown in Theorem 1 and in Examples 1 and 2. Examples are given of countable systems of oscillators with a rational commensurate collection of frequencies such that each of these systems does not have a periodic trajectory.  This effect can only occur in an infinite-dimensional phase space.

Linear flow trajectories on a finite-dimensional torus can be divided into two types. Every trajectory of a linear flow on a finite-dimensional torus is either periodic or its projection onto a two-dimensional symplectic subspace, is dense in the projection of the torus onto this subspace \cite{AA, Che}. We show that in the infinite-dimensional phase space there are at least three types of linear flow trajectories. The linear flow in Example 1 has neither a periodic trajectory nor a dense projection onto any two-dimensional symplectic subspace. Consequently, a trajectory of a linear flow on an infinite-dimensional torus can belong to at least one of three types. Namely, the periodic trajectories belong to the first type.
Trajectories which projection on some symplectic subspace are transitive in the projection of the invariant torus on this subspace form the second type.
The third type of trajectories includes non-periodic trajectories such that their projections on any symplectic subspace are not transitive.
Thus, it is shown that in the infinite-dimensional case there exists a new type of linear flow trajectories.

\section{Infinite-dimensional Hamiltonian system}

In the present paper the following description of infinite-dimensional Hamiltonian system is considered (\cite{Chernoff, Khrennikov, KS-18}).
{An}
infinite-dimensional Hamiltonian system is {a} triplet $(E,\omega ,h)$.
Here $E$ is an infinite-dimensional real separable Hilbert space, $\omega$ is a shift-invariant symplectic form (i.e. $\omega$ is a non-{degenerate} skew-symmetric bilinear form on the space $E$), and $h$ is a densely defined real-valued function which is {Frechet}-differentiable on {on a dense linear sub--manifolds
of the space $E$ .
}

A symplectic form $\omega$ in a real separable Hilbert space $E$ is called {\it native}
if there is an orthonormal basis (ONB)  $\{ e_k\}$ in the space $E$ such that
${\omega (e_{2k-1},e_{2j})=\delta _{j,k},\ k,j\in \mathbb N}$,
where $\delta _{j,i}$ is Kronecker symbol.
A basis with the above property is called symplectic basis of the symplectic space $(E,\omega )$.

A native symplectic form $\omega$ defines the decomposition of the space $E$ into orthogonal sum of two subspaces $E=P\oplus Q$ such that orthonormal systems $g_j=e_{2j-1},\, j\in \mathbb N$ and $f_k=e_{2k},\ k\in \mathbb N$ are ONB in subspaces $P$ and $Q$ respectively. Then
\begin{equation}\label{simp}
\omega (g_j,g_i)=0,\, \omega (f_i,f_j)=0\ \forall \ i,j\in {\mathbb N};
\quad \omega (g_j,f_k)=\delta _{jk},\, j,k\in \mathbb N .
\end{equation}
{In particular}, $\{ e_i,\ i\in \mathbb N \}=\{ g_j,f_k;\ j,\, k\in {\mathbb N}\}$ is {a} symplectic basis of the symplectic space $(E,\omega )$
(see \cite{Khrennikov, KS-18}).

The bilinear form $\omega $ on the Hilbert space $E$ is bounded. Therefore, it defines
{a} bounded linear operator ${\bf J}$ {associated} with the form $\omega$. {The} operator $J$ is {a} non-{degenerate}
skew-symmetric linear operator in the Hilbert space $E$.
{Conversely, any orthogonal basis $\{ g_j,f_k;\ j\in{\bf N},\, k\in {\bf N}\}$ of
the Hilbert space $E$ uniquely defines a skew-symmetric operator $\bf J$ on the Hilbert space $E$
by the equalities ${\bf J}(g_j)=-f_j,\ {\bf J}(f_k)=g_k,\  j\in {\bf N},\, k\in {\bf N}$.
The symplectic form $\omega $, associated to $J$, is native on the space $E$ with a symplectic basis $\{ g_j,f_k;\ j\in {\bf N},\, k\in {\bf N}\}$.
}

{The} subspaces  $Q$ and $P$ are called configuration space and momentum space respectively. The space $P$ has the role of the space which is conjugated to the space $Q$
(\cite{Kuksin, Khrennikov, KS-18, SSh-20}).

In a Hamiltonian system $(E,{\omega },h),$ the Hamilton function is {a} functional $h:\ E_1\to \mathbb R$ such that $E_1$ is a dense linear manifold of the Hilbert space $E$ and, moreover, the functional $h$ is {Frechet}-differentiable on {a}
linear manifold $E_2$ which is dense in the space $E$. The article {\cite{Kuksin}}
describes the properties of an infinite-dimensional Hamiltonian system and its invariant manifolds.

{The} Hamilton equation for {the} Hamiltonian system
$(E,{\omega },h)$ is the following equation $z'(t)=J(h'(z(t))),\, t\in \Delta,$ {in}
the unknown function $z:\ \Delta \to E_2$ where $\Delta $ is some interval of real line.  (\cite{Khrennikov, KS-18}).

A densely defined vector field ${\bf v}:\ E_2\to E$ is called Hamiltonian if there is
{a function $h:\ E_1\to \mathbb R$, {Frechet}-differentiable on the linear manifold $E_2$ and satisfying} the condition
${\bf v}(z)={\bf J}Dh(z),\, z\in E_2$.
Here $Dh$ is the differential of the function $h$ {and} $\bf J$ is the symplectic operator associated with the symplectic form $\omega $ in the space $E$.

A one-parameter group $g^t,\, t\in {\mathbb R},$ of continuously differentiable self-mappings of the space $E_2$ is called {a} smooth Hamiltonian flow in the space $E_2$ along the Hamiltonian vector field  ${\bf v}:\ E_2\to E$, if  ${{d }\over {d t}}g^t(q,p)={\bf v}(g^t(q,p)),\, (p,q)\in E_2,\ t\in \mathbb R$. If the Hamiltonian flow in the space  $E_2$ has {a}
unique continuous extension to the space $E$ then this extension of the flow is called generalized
Hamiltonian flow in the space $E$ along the Hamiltonian vector field
$\bf v$ (also this flow is called the generalized Hamiltonian flow generated by the Hamilton function $h$).

A symplectic operator $\bf J$ on the Hilbert space  $E$ defines {a} bijective
mapping ({called complexification}) of the space  $E$ to the complex Hilbert space
$H$ such that the multiplication on the imaginary unit $i$ in the space $H$ corresponds to the action of {the} operator $\bf J$ in the space $E$ (\cite{Khrennikov}). The inverse mapping $H\to E$ is called realization of a complex Hilbert space $H$ to the real Hilbert space $E$.
{A} realization of  a complex Hilbert space transforms a Schrodinger equation into the linear Hamilton equation (\cite{GS22, Khrennikov, KS-18, SSh-20}).

Conditions on the vector field $\bf v$ which {are} sufficient {for}
the existence of a flow in the space $E$ along the vector field  are considered in the monograph
S.G. Kreyn \cite{SGK}.
In the case of a Hilbert space $E$ and a linear vector field {a} sufficient condition {for the existence} of the flow is the self-adjointness of the linear operator
$z\to Dh(z)$ in the complexification of the space  $(E,{\bf J})$.
The sufficient conditions for the existence of a Hamiltonian flow which is generated in the space  $E_2$ by {a} nonlinear Schrodinger equation
are described in \cite{Bourgain, DF, Ponce, Zhidkov}.

\section {Representation of Hamiltonian dynamics in the form of dynamics of a system of oscillators}

According to \cite{Volovich19}, if a Hamiltonian system admits an invariant measure, then this system is {completely integrable in the Liouville sense}.
More precisely, there is an isomorphism of {the} above Hamiltonian system in the Koopman representation {with} the Hamiltonian system of a family of harmonic
oscillators.

{According to the definition in \cite{Volovich19}}, a system of harmonic oscillators (SHO) is a triplet
\begin{equation}\label{eq1}
(X,\mu, {\bf V}):
\end{equation}
where $(X,\mu )$ is a measurable space with a measure; $\bf V$ is {a $1$--parameter, strongly continuous sub--group of the}
special unitary group in the Hilbert space $L_2(X,\mu ,{\mathbb C})\equiv \cal H$,
${\bf V}_t: L_2(X,\mu )\to L_2(X,\mu ),\quad t\in {\mathbb R},$
{is defined} by the equality
\begin{equation}\label{2}
{\bf V}_t\phi (x)=e^{itf (x)}\phi (x),\ t\in \mathbb R .
\end{equation}
Here $f :\ X\to R$ is a measurable real-valued function, {and}
$\phi \in L_2(X,\mu )$.

{The} SHO (\ref{eq1}) is a Hamiltonian system $(E,\omega ,{h})$, where $E$ is the realization of {the}  complex Hilbert space $\cal H$ { with real
structure given by ${\cal H}= {\hbox{Re}}({\cal H}) \oplus i {\hbox{Im}}({\cal H})$, where $\oplus$ denotes
direct vector space sum,} which is defined by the equality
\begin{equation}\label{ove}
{\cal R}(\phi )=(q,p)=({\rm Re }\phi ,{\rm Im}\phi ),\ \phi \in \cal H.
\end{equation}
Here the symplectic operator $\bf J$ is induced by the complex structure of the space $H$ by the rule
${\bf J}({\cal R}(\phi ))={\cal R}(i\phi )$, therefore,  $\omega (\phi ,\psi)=({\bf J}{\cal R}(\phi ),{\cal R}(\psi ))_E={i\over 2}[(\phi ,\psi )_{\cal H}-(\psi ,\phi )_{\cal H}]$. The Hamilton function  $h:\ E\supset E_1\to \mathbb R$ is connected with the function $f$ {in the definition of SHO} by the equality
\begin{equation}\label{H}
h(p,q)={1\over 2}\int\limits_{X}f (x)(p(x)^2+q(x)^2)d\mu (x),\quad (p,q)\in E.
\end{equation}

The space ${\cal H}=L_2(X,\mu ,{\mathbb C})$ (and its realization) is called {the}
phase space of {the SHO}.
The unitary group $\bf V$ in the space $\cal H$ is {represented} in the realization space $E$ by the phase flow
\begin{equation}\label{P}
\Phi:\ \Phi _t={\cal R}{\bf V}_t{\cal R}^{-1}.
\end{equation}
The phase flow $\Phi $ is {a} group of orthogonal mappings such that every two-dimensional symplectic subspace of symplectic space $E$ is flow invariant.
In coordinates relative to the symplectic basis, the flow $\Phi$ is defined by the equality
$$
\Phi _t(q(x),p(x))=
$$
$$
=(q(x)\cos (tf (x))+p(x)\sin (tf (x)),p(x)\cos (tf (x))-q(x)\sin (tf (x))).
$$

{The} SHO (\ref{eq1}) and the initial point $\phi \in \cal H$ define a trajectory $\Gamma _{\phi}$ in phase space $\cal H$
\begin{equation}\label{3}
{\Gamma }_{\phi}=\{ e^{itf (x)}\phi (x),\ t\in {\mathbb R}\}.
\end{equation}

{The} problems of conservation l{a}ws and integrability of
linear systems in Hilbert space with additional positive integrals are considered in \cite{TSh, K22}.

In the case of SHO we have the following situation.
{The} Hamiltonian system (\ref{H}) has {a} collection of first
integrals $I_x:\ E\to {\mathbb R}$ {given} by the equalities
$$I_x(q,p)={1\over 2}f(x)(p(x)^2+q(x)^2),\quad x\in X.$$
{This} system of first integrals defines invariant manifolds of the Hamiltonian flow (\ref{P}) that can be parametrized by the collection of positive (nonnegative) $\mu $-measurable functions
$\{ r:\ X\to (0,+\infty )\}$
$$
{\mathbb T}_r=\{ I_x(q,p)=r(x),\ x\in X\}.
$$
{The} manifold ${\mathbb T}_r$ is called {an} invariant torus of {the} Hamiltonian system (\ref{H}).

For every nonnegative $\mu $-measurable function $r\in L_2(X,\mu )$ the invariant torus ${\mathbb T}_r$ of the SHO (\ref{eq1}) is the following set
\begin{equation}\label{tor}
{\mathbb T}_r=\{ u\in L_2(X,\mu ): \ |u(x)|=r(x),\, x\in X\}
\end{equation}
in the space $L_2(X,\mu)$.

\bigskip

Systems of harmonic oscillators and {the} corresponding invariant tori in a Hilbert space arise in consideration of quantum systems. A quantum system is a pair $({\cal H}, {\bf U})$ where  $\cal H$ is a separable complex Hilbert space and  $\bf U$ is a {$1$--parameter, strongly continuous} unitary group in the space $\cal H$.
The unitary group $\bf U $ of a quantum system has {a unique} generator $\bf L$ which is a self-adjoint operator in the space $\cal H$. Then, the unitary group is given by the equality ${\bf U}(t)=e^{-it{\bf L}},\, t\in {\mathbb R}$. The dynamics of this quantum system is described by the Schrodinger equation $i\dot \psi ={\bf L}\psi$.

In the case of a finite-dimensional Hilbert space, the equivalence of the Schrodinger equation to the equation of dynamics of a system of oscillators is described in \cite{MGI}.
For the infinite-dimensional case, using the spectral theorem, we obtain the following statement.

{\bf Theorem} \cite{Volovich19}.
{\it For any quantum system $({\cal H},{\bf U})$ there are a unitary mapping  ${\bf W}: \ {\cal H}\to L_2(X,\mu )$ and a measurable function $f:\ X\to \mathbb R$ such that the Schrodinger equation $i\dot \psi ={\bf L}\psi $ is unitary equivalent to the {dynamics} of the system of harmonic oscillators  $i{{\partial}\over {\partial t}}\phi (t,x)=f (x)\phi (t,x),\ x\in X,\ t\in {\mathbb R}$. This unitary equivalence {is given} by the equality ${\bf W}{\bf U}(t)={\bf V}_t{\bf W},\ t\in {\mathbb R}$. It follows that any quantum system is completely integrable and an infinite number of integrals of motion are indicated.}

\section{Dynamics on  an infinite-dimensional torus}\label{Dyn-on-inf-dim-tori}

{In this paper we consider the dynamics of an SHO in the special case such that the measure $\mu$ in (\ref{eq1}) is a discrete measure. }
{We study SHOs satisfying the following conditions: $X=\mathbb N$ and $\mu $ is the counting measure on the set $\mathbb N$}. In this case ${\cal H}=L_2(X,\mu, {\mathbb C})=l_2$. A non-{degenerate} invariant torus (\ref{tor}) in this case is defined by the sequence of
positive numbers $\{ r_n\}\in l_2$. This torus is the following set
\begin{equation}\label{tn}
{\mathbb T}=\{ u\in l_2:\ |u_n|=r_n,\, n\in {\mathbb N}\}.
\end{equation}

A realization ${\cal R}$ (\ref{ove}) of the complex Hilbert space $\cal H$ mappings a point
$u\in L_2({\mathbb N},\mu, {\mathbb C}))$ to the point  $(q,p)\in E=Q\oplus P$, where  $Q=P=L_2({\mathbb N},\mu ,{\mathbb R})$, and the following system of equalities are satisfied
$u_n=q_n+ip_n \ \forall \ n\in {\mathbb N}$. The realization of the torus  (\ref{tn}) is the following set in the space $E$:
\begin{equation}\label{tnr}
{\mathbb T}=\{ (q,p)\in E:\ p_n^2+q_n^2=r_n^2,\, n\in {\mathbb N}\}.
\end{equation}

An infinite-dimensional torus (\ref{tnr}) in the Hilbert space $E=Q\oplus P$ is the countable Cartesian product of circles
\begin{equation}\label{c2}
\mathbb T=C_1\times C_2\times.... \subset E.
\end{equation}
Here $C_k=\{ p_k^2+q_k^2=r_k^2\} $ is the circle in the two-dimensional {plane}
$E_k={\rm span }(f_k,g_k)$. Since $\mathbb T\subset E$,
\begin{equation}\label{32}
\sum\limits_{k=1}^{\infty }r_k^2<+\infty .
\end{equation}

The Hamilton function
\begin{equation}\label{HNI}
H(q,p)=\sum\limits_{k=1}^{\infty }{{\lambda _k}\over 2}(p_k^2+q_k^2), \quad (q,p)\in E_1,
\end{equation}
where $\{ \lambda _k\}$ be a sequence of positive numbers and $E_1=\{ (q,p)\in E:\ \sum\limits_{k=1}^{\infty }\lambda _k(q_k^2+p_k^2)<+\infty \}$, defines
{a} flow on the torus (\ref{tnr}) in the space
$E$ by the following equalities
\begin{equation}\label{Z}
z_k(t)=({\bf V}_tz_{0})_k=e^{i\lambda _kt}z_{k0},\ t\in \mathbb R;\quad z_k=q_k+ip_k.
\end{equation}

This Hamiltonian flow should be considered firstly in the subspace
$E_2=\{ x\in E:\ \sum\limits_{k=1}^{\infty}\lambda _k^2(q_k^2+p_k^2)<\infty \}$ of the space  $E=l_2$. The reason is the following. The tangent vectors of the flow is defined by the gradient of the Hamilton function (\ref{HNI}) and this function is differentiable  on the subspace $E_2$ in its domain $E_1$.

Let us study {the} existence of periodic trajectories
$\Gamma _{z_0}=\{ z(t),\ t\in \mathbb R \}$, of the Hamiltonian system (\ref{HNI}) on the torus $\mathbb T$.

{\bf Lemma 1}.
{\it Let  $\Gamma $ be a periodic trajectory of the Hamiltonian system (\ref{HNI}) on the torus $\mathbb T$ with a period  $\tau$. Then, for every  $n\in \mathbb N$ the projection  $\Gamma _n$  of the trajectory $\Gamma $ onto the $2n$-dimensional subspace $\oplus_{k=1}^nE_k$ is the periodic trajectory on the torus
${\mathbb T}_n$. Here ${\mathbb T}_n$ is the projection of the torus $\mathbb T$ onto  $n$-dimensional subspace $\oplus_{k=1}^nE_k$. Moreover, the period of {the} projection of a trajectory {is not} greater than the period of trajectory.}

{\bf Lemma 2.}
{\it Let the sequence $\{ {{2\pi }\over{\lambda _k}}\} $
of periods of oscillators in the Hamiltonian system (\ref{HNI}) be unbounded. Then
{any} trajectory of this Hamiltonian system
on the non-{degenerate} torus (\ref{tnr})
is not periodic.
 }

\begin{proof}

Let us assume the opposite. Then there is a trajectory $\Gamma$ of
{the} Hamiltonian system (\ref{HNI}) which is periodic with {period}
$\tau \in (0,+\infty )$. A trajectory $\Gamma$ belongs to some torus ${\mathbb T}\subset l_2$. Then, according to {Lemma} 1, for every $m\in\mathbb N${,} the projection $\Gamma _m$ onto {the} two-dimensional subspace  ${\rm span}(e_{2m-1},e_{2m})$ is
{a} periodic trajectory on the torus ${\mathbb T}_m$ which is the projection of the torus  $\mathbb T$ onto  two-dimensional subspace ${\rm span}(e_{2m-1},e_{2m})$. Moreover, the period of trajectory $\Gamma _m$ is no greater than  $\tau$.
Hence, $\tau \geq \tau _m\ \forall \ m$.  Since the sequence $\{ \tau _m\}$ is unbounded, there is no period $\tau $ of the trajectory $\Gamma$.
\end{proof}

{\bf Lemma 3.}
{\it {Suppose that} the sequence $\{ \lambda _k\}$ satisfies the rational commensurability condition, i.e. there is a finite collection of numbers
$\lambda _{j_1},...,\lambda _{j_m}$ such that one of them is a linear combination,
{with rational coefficients, of the others}.
Then any trajectory of Hamiltonian system (\ref{HNI}) is not dense everywhere in a non-{degenerate} torus $\mathbb T$. }

\begin{proof}
Let $\Gamma$ be a trajectory of the Hamiltonian system (\ref{HNI}) which belongs to a non-{degenerate} torus $\mathbb T$. According to the Jacobi theorem and to the assumptions of Lemma 3,
the projection
 $\Gamma _m$ of the trajectory $\Gamma$ onto $2m$-dimensional subspace ${\rm span}(e_{j_1},...,e_{j_{2m}})$ is not dense everywhere in the projection ${\mathbb T}_m$ of the torus  $\mathbb T$. Therefore, the trajectory  $\Gamma$ is not dense everywhere in the torus $\mathbb T$.
\end{proof}

The following theorem is {a} consequence of {Lemmas} 2, 3.

{\bf Theorem 1.}
{\it Let the sequence $\{ \lambda _k\}$ is rationally commensurable, i.e. there is a finite
collection of numbers $\lambda _{j_1},...,\lambda _{j_m}$ such that one of them is a linear combination, {with rational coefficients, of the others}.
Let the sequence $\{ {1\over {\lambda _k}} \}$ {be} unbounded. Then for any non-{degenerate} torus $\mathbb T$ every trajectory of Hamiltonian system (\ref{HNI}) is neither periodic, nor dense in the torus $\mathbb T$. }

{\bf Lemma 4.}
{\it Let Hamiltonian system in symplectic space  $(E,\ \omega )$ be defined by the Hamilton function (\ref{HNI}). Let ${\mathbb T}_r$ be a torus in the space $\cal H$ which is defined by the sequence of positive numbers $\{ r_k\} \in l_2$. If there is a periodic trajectory of the Hamiltonian system (\ref{HNI}) in the torus  ${\mathbb T}_r$, then every finite collection of numbers $\{ \lambda _{j_1},...,\lambda _{j_m}\}$ is strongly rationally commensurable.}

The statement of Lemma 4
is the consequence of Lemma 1 and the Jacobi theorem.

{\bf Lemma 5.}
{\it Let {the} Hamiltonian system in the symplectic space  $(E,\ \omega)$
{be} defined by the Hamilton function  (\ref{HNI}).
Let ${\mathbb T}_r$ be a torus in the space $\cal H$ {defined} by the sequence of positive numbers $\{ r_k\} \in l_2$. If there is a periodic trajectory on the torus ${\mathbb T}_r$ then the sequence of frequencies $\{ \lambda _k\}$ has the form $\{ \lambda _0{n_k}\}$ where  $\lambda _0>0$ and $\{{n_k}\}$ is {a} sequence of natural numbers. }

\begin{proof} Let us note by $ T_n$ the minimal period of the projection of a trajectory  {onto the}
symplectic subspaces $\oplus_{k=1}^nE_k$ for every $n\in \mathbb N$. 
Then $T_1={{2\pi }\over {\lambda _1}}$ is the minimal period of the first oscillator in the phase plane $E_1$. Then the sequence $\{ T_n\}$ 
 is nondecreasing and the period
$T_n$ is {a} multiple of the period $T_1$ (since $T_n\lambda _1=2\pi m_n$ with some $m_n\in \mathbb N$) for every  $n\in \mathbb N$.

Moreover, \begin{equation}\label{karat}\forall \ n\in {\mathbb N}\ \exists \ l\in {\mathbb N}:\ m_{n+1}=lm_n.
\end{equation} 
Because, $T_{n+1}=lT_{n}$ since $T_n$ is the least period and $T_{n+1}$ is the period of the projection of a trajectory on the torus ${\mathbb T}_n$.

According to the Theorem 1 
the boundedness of the sequence  $\{ T_n\}$ is {a} necessary condition {for} the periodicity of a trajectory on the torus ${\mathbb T}_r$. Hence,
the sequence $\{ T_n\}$ is bounded. Thus, it is {a} stationary sequence {(i.e. its elements has a constant values for sufficiently large numbers)}.
Since $T_1$ is the minimal value of the monotone sequence $\{ T_n\}$, the condition
$T_n=m_nT_1\ \forall \ n\in \mathbb N$ holds where  $\{ m_n\}$ is nondecreasing bounded sequence of natural numbers. Hence, there is a natural number  $M=\max \{ m_k\}$ such that $T_n=MT_1$ for every sufficiently large $n$. As we already have seen, if $k\in \mathbb N$ then $T_n$ is the multiply of ${{2\pi }\over {\lambda _k}}$ for every $n\geq k$. Therefore, there is a natural number $n_k$ such that
${{2\pi }\over {\lambda _k}}={{MT_1}\over {n_k}}$. Therefore, $\lambda _k=\lambda _0n_k$ where $\lambda _0={{2\pi }\over {MT _1}}.$
\end{proof}

{\bf Corollary 1.}  {\it If  there is a periodic trajectory $\Gamma$ of the Hamiltonian system (\ref{HNI}) has the period $T$, then $T $
is a multiple of the period  $T _m$ of the projection of the trajectory $\Gamma _m$ onto the symplectic subspace ${\rm span}(e_1,e_2,...,e_{2m})$ for every $m\in\mathbb N$.}

As it is shown in the proof of Lemma 5, $T=T_ml_m$ with some $l_m\in \mathbb N$ according to (\ref{karat}).

{\bf Theorem 2.}
{\it Let $\{ \lambda _k\}$ {be} the sequence of frequencies of the Hamiltonian system
(\ref{HNI}). Let $r\in l_2$ be the sequence of positive numbers.
Then every trajectory of the Hamiltonian system (\ref{HNI}) on the torus ${\mathbb T}_r$ is periodic if and only if the sequence  $\{ \lambda _k\}$ satisfies the following condition.
There are the number
$\lambda _0>0$ and the  sequence of natural numbers
$\{ n_k\}$ such that $\lambda _k=\lambda _0n_k,\ k\in \mathbb N$. }

\begin{proof}
According to Lemmas 4, 5
the condition of Theorem 2
is necessary {for} the periodicity of  a trajectory of the Hamiltonian system (\ref{HNI}).

Let the condition of the Theorem 2
{be satisfied}. Let $T={{2\pi }\over {\lambda _0}}$.
Then the number $T$ is the period of the projection of a trajectory onto a two-dimensional symplectic subspace  $E_k={\rm span}(f_k,g_k)$ for every $k\in \mathbb N$. Therefore, any trajectory of SHO has {period} $T$.
\end{proof}

{\bf Remark}.
Theorem 2
presents the criterion for the periodicity of a trajectory on an infinite-dimensional (countable-dimensional) torus. This criterion is the analogue of the first part of the Jacobi theorem. The condition of Theorem 2
on the sequence of frequencies of a countable SHO is coincides with the condition of strong rational commensurability in the Jacobi theorem in the case of finite set of values of the sequence. In the case of infinite set of values the condition of {Theorem} 2
implies the condition of strong rational commensurability.
Theorem 1 shows that for the periodicity of a trajectory, only one condition of strong rational commensurability is not enough.

{\bf Example 1}. If $\lambda _k={1\over {k!}},\ k\in \mathbb N$, then the countable family of frequencies $\{ \lambda _k,\ k\in \mathbb N\}$ {satisfies} the strong rational commensurability condition. But according to Lemma 2
there is no periodic trajectory of the SHO.

It may seem that this effect is associated with the tendency of frequencies to zero. But this is not true, as the following example shows.

{\bf Example 2}. Let us show another example of rationally commensurable family of frequencies which are {different} from zero and infinity {and} such that there is no
periodic point of {the} SHO.
Let $\{ p_k\}$ be a sequence of prime numbers numbered in ascending order. Then the sequence of frequencies with the needed properties can be chosen in the following form  $\lambda _k={{p_{2k}}\over {p_{2k-1}}},\ k\in \mathbb N$. {A trajectory of this SHO is not periodic since periods of trajectories of its finite subsystems of harmonic oscillators are unbounded. In fact, the period of the oscillator with number $k$ is equal to $\tau _k={{2\pi p_{2k-1}}\over {p_{2k}}}$ and the period of the system of first $N$ oscillators is equal to $\tau _{1,2,...,N}=2\pi p_1p_3...p_{2N-1}$.}

{In both cases of SHO in examples 1 and 2, the trajectory of the linear flow is not transitive in the corresponding invariant torus $\mathbb T$. Moreover, according to the first part of Jacobi's theorem, for any finite-dimensional symplectic subspace $E'$ the projection of a trajectory onto the subspace $E'$ is a periodic trajectory in the torus ${\mathbb T}'$ which is the projection of the torus $\mathbb T$ onto the subspace $E' $. According to Jacobi's theorem, in the finite-dimensional case, each trajectory of a linear flow on the corresponding invariant torus is either periodic or has a projection onto some finite-dimensional symplectic subspace, dense in the projection of the invariant torus onto this subspace. Thus, we observe a new type of linear flow trajectories on an infinite-dimensional torus. Trajectories of this type are neither periodic nor have projections onto some symplectic subspace that is dense in the projection of the invariant torus onto this subspace.}

\bigskip

Let us study the existence of a trajectory $\Gamma _{z_0}=\{ z(t),\ t\in \mathbb R \}$, of the Hamiltonian system (\ref{HNI}) such that this trajectory is dense everywhere in a torus
$\mathbb T$.

{\bf Lemma 6.}
{\it A trajectory $\Gamma $ is dense everywhere on a non-{degenerate} torus $\mathbb T$ if and only if, {for every $n\in \mathbb N$, the} projection  $\Gamma _n$ of the trajectory  $\Gamma $ onto subspace $ E_{(2n)}= {\rm span}(f_1,g_1,...,f_n,g_n)$ is dense
everywhere on the torus ${\mathbb T}_n$. }

\begin{proof}
The distance between {the} projections of two points of phase space $E$ onto the subspace $E_{(2n)}$ is {not} greater than the distance between the pre--image of this projections.
Therefore, the density everywhere on the torus ${\mathbb T}_n$ of the projection $\Gamma _n$ for some $n\in \mathbb N$  is {a} necessary condition for a density everywhere on the torus
${\mathbb T}$ of the trajectory $\Gamma$.

To prove the sufficiency of this condition{,} let us fix a number $\epsilon >0$.
Let us show that, {for any  $z\in {\mathbb T}$,}
$B_{\epsilon}(z)\bigcap \Gamma \neq\emptyset$,
where $B_{\epsilon }(z)=\{ u\in E:\ \|u-z\|_E<\epsilon \}$.
According to {condition} (\ref{32}) there is a number $N$ such that $\sum\limits_{k=N+1}^{\infty }r_k^2<{1\over {4}}\epsilon ^2$. By {assumption}
the projection
$\Gamma _N$ is dense everywhere in the torus  ${\mathbb T}_N$. By the choice of the number $N$ the distance between a point $z\in \mathbb T$ and its projection onto the torus
${\mathbb T} _N$ is no greater than
$\left( \sum\limits_{k=N+1}^{\infty }|z_k|^2\right)^{1/2}<{{\epsilon }\over 2}$. Therefore, for every $z\in \mathbb T$ there is a point $u\in \Gamma$ such that $\|u-z\|_E\leq \epsilon$.
\end{proof}

{\bf Theorem $3$}.
{\it A trajectory of {the} Hamiltonian system
(\ref{HNI}) is dense everywhere on the torus  $\mathbb T$ in the Hilbert space
$E$
if and only if the set of numbers $\{ \lambda _k,
, k\in {\mathbb N}\}$, is not rationally commensurable.}

\begin{proof}
The collection of numbers $\{ \lambda _k,\ k\in {\mathbb N}\}$ is not rationally commensurable if and only if for every $N\in {\mathbb N}$ the collection of numbers $\{ \lambda _1,...,\lambda _N\}$ satisfies the condition
\begin{equation}\label{10p}
\{\sum\limits_{k=1}^N\omega _{k}n_k=0,\quad n_k\in {\mathbb Z}\} \quad
\Rightarrow \quad \{ n_k=0,\ j=1,...,N\}.
\end{equation}

The necessity of the condition (\ref{10p}) is the consequence of Lemma 6
and second part of the Jacobi theorem.

Let us prove {sufficiency}.
Let $\Gamma $ be a trajectory of the Hamiltonian system (\ref{HNI}). Let $x_0\in \mathbb T$.
Let $U$ be a neighborhood of the point $x_0$. Then, there is a number $\epsilon >0$ such that the ball $B_{\epsilon }(x_0)$ belongs to the neighborhood $U$.
According to the condition ({\ref {32}}) there is a natural number $N$ such that $\sum\limits_{k=N+1}^{\infty}r_k^2<{1\over 4}\epsilon ^2$. The Jacobi theorem states that according to the condition (\ref{10p}) the projection of the trajectory $\Gamma$ to the subspace ${\mathbb R}^{2N}$ intersects the ball $B_{\epsilon /2}(x_0)$. Therefore, the trajectory  $\Gamma $ intersects the neighborhood $U$ of the point $x_0$. Thus, the trajectory $\Gamma $ is transitive in the torus $\mathbb T$.
\end{proof}

{\bf Remark.}
Theorem 3
is {a} criterion of density everywhere in {an} invariant torus  of
a countable SHO with {Hamilton function} (\ref{HNI}).
This criterion is the {analogue} of the second part of {Jacobi's}
theorem.

\section{Ergodicity of Kolmogorov's measure with respect to the Hamiltonian flow}

{Let us study the ergodic property of Kolmogorov's measure on the
infinite-dimensional torus} with respect to linear flow (\ref{Z}). Let us note that this measure is used in harmonic {analysis} on infinite-dimensional {tori}
\cite{Fufaev, Platonov}.

A torus $\mathbb T$ is {a} countable Cartesian product (\ref{tor}) of circles $C_k$ with radii $r_k$ such that {condition} (\ref{32}) holds. Let us equip every circle $C_k$ with the angular coordinate $\phi _k\in [0,2\pi ]$ and with the probability measure $\lambda _k$ which is multiple to the Lebesgue measure on the circle  $C_k$.
Then, the torus  (\ref{tor}) can be considered as the countable product of probability spaces $(C_k, {\cal B}(C_k), \lambda _k)$. Thus, the torus $\mathbb T$ is equipped with the probability measure $\lambda _{\mathbb T}=\otimes _{k=1}^{\infty }\lambda _k$. According to Kolmogorov's theorem on the product of probability spaces (see \cite{Bogachev}), the measure $\lambda _{\mathbb T}$ is defined on the $\sigma $-algebra {${\cal A}$ which} is generated by the countable product of  $\sigma$-algebras ${\cal B}(C_k),\ k\in \mathbb N$ of Borel subsets of circles $C_k$ and this measure $\lambda _{\mathbb T}$  is {a} probability measure on the torus $\mathbb T$. The $\sigma $-algebra ${\cal A}$ is generated by the collection of cylindrical subsets of the torus $\mathbb T$. The measure $\lambda _{
\mathbb T}:\ {\cal A}\to [0,1]$ is called Kolmogorov's measure.

Every circle $C_k$ is parametrized by the angular parameter  $\phi _k\in [0,2\pi ]$. Then, the torus $\mathbb T$ is parametrized by the vector valued parameter $\phi =\{ \phi _k\}\in \otimes_{k=1}^{\infty}[0,2\pi ]$.

A function $f:\ T\to \mathbb R$ is integrable with respect to the measure $\lambda _{\mathbb T}$ if and only if {it} can be approximated by {linear combinations}
of cylindrical {functions}.
Let
\begin{equation}\label{systy}
f(\phi )=\otimes_{k=1}^{\infty }u_k(\phi _k),\ \phi =\{ \phi _k\}\in \otimes_{k=1}^{\infty}[0,2\pi ],
\end{equation}
be a tensor product of functions  $u_k\in L_2({[0,2\pi ]})$, $k\in \mathbb N$. The condition of approximation of a function $f$ by {a} linear combination of cylindrical
{functions} implies that the tensor multipliers $u_k$ {tend}
to the unit function.

The space $H({\mathbb T})=L_2({\mathbb T},{\cal A},\lambda _{\mathbb T},{\mathbb C})$ is defined (\cite{Fufaev}) as the completion of linear {span} of {the}
system (\ref{systy})  by the Euclidean norm
\begin{equation}\label{scala}
\| \otimes_{k=1}^{\infty }u_k(\phi _k)\|_{H({\mathbb T})}=\prod\limits_{k=1}^{\infty }\| u_k\|_{L_2(C_k,{\cal B}(C_k),\lambda _k,{\mathbb C})}
\end{equation}

The function
\begin{equation}\label{FF}
u(\phi )=\exp(i({\bf \phi},{\bf n}))=\exp( i \sum\limits_{k=1}^{\infty }\phi_kn_k),\ x\in T,\quad {\bf n}\in c_0({\mathbb N},{\mathbb Z}),
\end{equation}
where ${\bf n}:\ {\mathbb N}\to {\mathbb Z}$, is correctly defined if and only if the sequence of integers $\bf n$ is finite.
If the sequence of integers $\bf n$ is finite then  (\ref{FF})
 is {a} cylindrical function.

The system of functions (\ref{FF}) forms {a} basis in the space
$H({\mathbb T})=L_2({\mathbb T},{\cal A},\lambda _{\mathbb T},{\mathbb C})$ (\cite{Fufaev}),
since it is orthogonal and complete according to the definition oh the space  $H({\mathbb T})$.

Recall that {a} probability measure $\lambda$ on the measurable space
$(\mathbb T ,{\cal A})$ is called ergodic with respect to {a} family $\cal V$ of measurable mappings of the
set $\mathbb T$ into itself  if there is no $\cal V$-invariant subset $A\in \cal A$ such that
$\mu (A)\in (0,1)$.  The ergodicity of the measure with respect to {a} family
${\cal V}$ of measurable {maps} of the set $\mathbb T$ into itself is equivalent to the following condition. If a measurable function $f: \mathbb T \to \mathbb R$ is invariant with respect to mappings
$f\to f\circ {\bf V}^{-1},\, {\bf V}\in \mathcal V,$ then the function $f$ is constant up to
{a} subset {of} zero $\lambda$-measure \cite{KSF}.

{\bf Theorem 4.}
{\it
The ergodicity of Kolmogorov's measure $\lambda _{\mathbb T}$
with respect to the flow (\ref{Z}) is equivalent to the condition of density everywhere on an invariant torus  $\mathbb T$ of a trajectory of the flow (\ref{Z}).
}

\begin{proof}
According to {Theorem} 3{,}
it is sufficient to prove that {condition} (\ref{10p}) of {everywhere density of a} trajectory is equivalent to the condition of ergodicity of Kolmogorov's measure on the torus  $\mathbb T$ with respect to the flow  (\ref{Z}).

1) Necessity. If there is a natural number $N$ such that the condition (\ref{10p}) is violated, then there is a function of the type (\ref{FF}), which is invariant with respect to the flow (\ref{Z}) but is not constant.

2) Sufficiency.
Let condition (\ref{10p}) {be} satisfied for arbitrary natural $N$.
If the function {$u\in H({\mathbb T})$ is} invariant with respect to the mappings  $u\to u\circ {\bf V}_t^{-1},\, t\in \mathbb R,$ then
\begin{equation}\label{Four}
(u\circ {\bf V}_t^{-1},u_{\bf n})_{H({\mathbb T})}=(u,u_{\bf n})_{H({\mathbb T})}
\end{equation}
for every finite sequence of integers $\bf n$.
Since the system of functions  (\ref{FF}) is {an} ONB in the space  $H({\mathbb T})$, then $u=\sum\limits_{j=1}^{\infty }\alpha _ju_{{\bf n}_j}$, where
$\{\alpha _j\}\in l_2$ и $\{ {\bf n}_j\}: \ {\mathbb N}\to c_0({\mathbb N},{\mathbb Z})$.
Therefore, the invariance condition
(\ref{Four}) implies that $\alpha _j(1-e^{it(\lambda ,{\bf n}_j)_{l_2}})=0\ \forall \ j\in {\mathbb N},\, t\in {\mathbb R}$. {This} is possible if and only if $(\lambda ,{\bf n}_j)_{l_2}=0$ for every nontrivial {component} in the decomposition of the function $u$. According to {condition} (\ref{10p}), this means that only the component with the trivial vector  ${\bf n}_j$ can be nonzero. Therefore, under
{assumption} (\ref{10p}) only {the almost everywhere constant
functions} $u\in H({\mathbb T})$ {are} invariant with respect to {the}
mappings $u\to u\circ {\bf V}_t^{-1},\, t\in \mathbb R$.
\end{proof}

{\bf Remark.}
Theorem 4
states that the condition of  density everywhere on an invariant torus of a trajectory of countable SHO with the Hamilton function (\ref{HNI}) is equivalent to the ergodicity of Kolmogorov's measure on the torus with respect to the Hamiltonian flow. Thus, the condition   (\ref{10p}) of absence of rational commensurability of the family of frequencies is the criterion of ergodic randomness of the dynamics of {the} SHO  (\ref{HNI}).

\section
{Nonwandering property}

Let us remember that a point $A$ of a trajectory  $\Gamma $ is called {a} nonwandering point with respect to the flow ${\bf V}_t,\ t\geq 0$ if for any neighborhood of the point $A$ and for any number $T>0$ there is a number $t$ such that  $t>T$  and $O(A)\bigcap {\bf V}_t(O(A))\neq \emptyset$ \cite{Efremova}.

{\bf Theorem 5.}
{\it Every point of the phase space $E$ is {a} nonwandering point of the flow of the Hamiltonian system (\ref{HNI}).
}

\begin{proof} Let $z\in E$. Then, there is {a} torus ${\mathbb T}_r$ such that
$z\in {\mathbb T}_r$. Let us fix a number $\epsilon >0$. Then, there is a number $m\in \mathbb N$ such that $\sum\limits_{k=m+1}^{\infty }r_k^2<{{\epsilon ^2}\over 8}$. Therefore, for every point of the torus ${\mathbb T}_r$ the distance between this point and its orthogonal projection on the subspace $F_m=\{ (p,q)\in E:\ p_i=0,\, q_i=0,\ \forall \ i=1,...,m\}$ is {not}
greater than  ${{\epsilon }\over 2}$.

According to Poincare's {recurrence} theorem for {a} dynamics in
{an} invariant finite dimensional subspace, the orthogonal projection
{$F_m^{\bot}$ of} $z_m$ of the point $z$ {into} the subspace $F_m^{\bot}$ {is a} nonwandering point of the flow (\ref{Z}) in the invariant torus
${\mathbb T}_m={\mathbb T}\bigcap F_m^{\bot}$.
Hence, the point $z$ is {a} nonwandering point of the flow (\ref{Z}).
\end{proof}

\section{Conclusion}
In present paper, {we generalize Jacobi's and Weyl's theorems to linear flows
on infinite-dimensional tori}.
As in the {above quoted} classical theorems on the phase flow of a finite system
of oscillators on a finite dimensional torus, conditions are obtained for periodicity, transitivity and ergodicity {and density of the trajectories of an infinite system of oscillators
on an invariant torus.
}
The property of ergodicity of the linear flow on an invariant torus of a countable system of harmonic oscillators and the property of nonwandering of a point of phase space are studied.
We prove the ergodicity of a measure (and the ergodicity of the space of functions) on an invariant torus of the flow of an infinite system of oscillators. The nonwandering property of an arbitrary point of the phase space of infinite-dimensional linear flow of a countable system of Harmonic oscillators are proven.
It is shown that for infinite-dimensional linear flows there is a new type of trajectory that is absent in the finite-dimensional case.

\section{Acknowledgments} Authors thank I.Ya Aref'eva and L.S. Efremova for the interest in this work and useful discussion. 
We especially thank Luigi Accardi for his careful reading of the manuscript and many important comments that helped improve and clarify the presentation.

\end{document}